\documentclass{amsart}
\usepackage{amsmath,amsthm, amscd, amssymb, amsfonts}

\numberwithin{equation}{section}
\begingroup
\newtheorem{propo}[equation]{Proposition}
\newtheorem{corol}[equation]{Corollary}
\newtheorem{theor}[equation]{Theorem}
\newtheorem{lemma}[equation]{Lemma}
\theoremstyle{definition}
\newtheorem{defin}[equation]{Definition}

\newcommand{\al}{\alpha}
\newcommand{\be}{\beta}
\newcommand{\gm}{\gamma}
\newcommand{\sg}{\sigma}
\newcommand{\io}{\iota}

\newcommand{\Dt}{\Delta}
\newcommand{\ze}{\zeta}
\newcommand{\la}{\lambda}
\newcommand{\ZZ}{\mathbb{Z}}
\newcommand{\FF}{\mathbb{F}}

\newcommand{\Aut}{\operatorname{Aut}}

\newcommand{\Gal}{\operatorname{Gal}}
\newcommand{\GL}{\operatorname{GL}}
\newcommand{\id}{\operatorname{id}}

\def\orb{{\mathcal O}}
\newcommand{\trid}{\triangleright}
\newcommand{\mdpd}[4]{\left(\begin{smallmatrix}#1 & #2 \\ #3 & #4\end{smallmatrix}\right)}
\newcommand{\zpzp}{\ZZ_p\oplus\ZZ_p}
\newcommand{\Ss}{\mathbb{S}}

\begin{document}
\title{Indecomposable racks of order $p^2$}
\author{Mat\'\i as Gra\~na}
\address{Mat\'\i as Gra\~na \newline
\indent MIT, Mathematics Department \newline
\indent 77 Mass. Ave. \newline
\indent 02139 Cambridge, MA - USA}
\email{matiasg@math.mit.edu}
\begin{abstract} 
We classify indecomposable racks of order $p^2$ ($p$ a prime). There are $2p^2-2p-2$
isomorphism classes, among which $2p^2-3p-1$ correspond to quandles. In particular,
we prove that an indecomposable quandle of order $p^2$ is affine. One of the
results yielding this classification is the computation of the quandle nonabelian
second cohomology group of an indecomposable quandle of prime order; which turns out
to be trivial, as in the abelian case.
\end{abstract}
\thanks{This work was supported by CONICET}
\maketitle

\section{Introduction and Notation}
Racks and quandles have been considered since 1982 by Joyce \cite{j}
and Matveev \cite{ma} as a tool to provide knot invariants. Since then,
cohomology theories for them have been introduced independently at least with
two main goals: to give new knot invariants and invariants of knotted surfaces
on one hand (cf. \cite{cjkls}), and to classify Yetter--Drinfeld modules over groups
on the other hand (cf. \cite{g}). These theories coincide, which put racks,
as a subject, between topology and Hopf algebras.
We refer to the survey \cite{os} for topological issues and to
\cite{ag} for Hopf algebraic ones.

Let us define the object of study:
\begin{defin}
A \emph{rack} is a pair $(X,\trid)$ where $X$ is a set and $\trid:X\times X\to X$ is a binary
operation satisfying
\begin{align}
\label{Q1}&\text{The functions }\phi_x:X\to X,\ 
	\phi_x(y)=x\trid y\text{ are bijections for all }x\in X, \\
\label{Q2}&x\trid(y\trid z)=(x\trid y)\trid(x\trid z)\ \forall x,y,z\in X, \\
\intertext{A rack is a \emph{quandle} if it further verifies}
\label{Q3}&x\trid x=x\ \forall x\in X.
\end{align}
The main model for racks are unions (possibly with repetitions) of conjugacy classes
in a group, where the operation is the conjugation $x\trid y=xyx^{-1}$. When there
are repetitions, $x\trid y$ is defined to be in the same copy as $y$. Notice that
any of these racks is actually a quandle.
\end{defin}

So far, the classification of (finite) racks turns out to be a natural problem,
though it seems to be out of reach. Classifications which seem more tractable
include that of racks of low order (which can be either done by computers), 
that of indecomposable racks, of faithful racks (see definitions below).
A first approach in this direction was done in \cite{egs},
where indecomposable set-theoretical solutions to the
Yang-Baxter equation with a prime number of elements were classified. As a
byproduct, it included the classification of indecomposable racks of prime order.

Here we classify indecomposable racks of order $p^2$ ($p$ a prime); the main
techniques are contained in \cite{egs} and \cite{ag}.
This is the natural place to put the \\
\textbf{Acknowledgments.} I thank N. Andruskiewitsch and P. Etingof for valuable comments.
I also thank T. Ohtsuki who, by giving a list of indecomposable quandles of order
$9$, encouraged me to write this paper. Finally, I thank the warm hospitality of MIT
and its productive atmosphere in which I wrote it.

\medskip
We give now the definition of indecomposability, as well as the necessary
definitions for reading the paper.
\begin{defin}
A rack $(X,\trid)$ is \emph{decomposable} if it can be split properly into stable subracks, i.e.,
if $X=Y\sqcup Z$ (a disjoint union), neither of them empty, and $X\trid Y=Y$, $X\trid Z=Z$.
A rack is \emph{indecomposable} if it is not decomposable.

Let $\phi:X\to\Ss_X$ be the function defined in \eqref{Q1}. We denote by $G_X^0$ the subgroup
of $\Ss_X$ generated by the image of $\phi$. This group operates on $X$ by rack automorphisms;
i.e, if $\sg\in G_X^0$ then $\sg(x\trid y)=\sg(x)\trid\sg(y)$. A rack is \emph{faithful} if
$\phi$ is injective. In this case $X$ is naturally seen as a union of conjugacy classes inside
$G_X^0$, and then it is a quandle.
If $X$ is a rack, we denote by $\Aut(X)$ the group of rack automorphisms of $X$.

An \emph{affine quandle} (also called \emph{Alexander quandle}) is a pair $(A,g)$,
where $A$ is an abelian group and $g\in\Aut(A)$.
It is a quandle with the structure $x\trid y=(1-g)(x)+g(y)$. It is easy to see that such a
quandle is indecomposable iff $1-g$ is surjective and it is faithful iff $1-g$ is injective.
If $A$ is cyclic, we denote $g$ usually by $g(1)$. It is proved in \cite{egs} that an indecomposable
quandle of prime order $p$ is affine, isomorphic to $(\ZZ_p,q)$, where $q\in\ZZ_p^*-\{1\}$.

If $X$ is a quandle and $H$ is a group, the \emph{non-abelian $2$-cocycles with values in $H$}
are the functions $\be:X\times X\to H$ such that
$$\be(x,y\trid z)\be(y,z)=\be(x\trid y,x\trid z)\be(x,z)\quad\forall x,y,z\in X.$$
Two cocycles $\be$ and $\be'$ are \emph{cohomologous} if there exists a function
$\gm:X\to H$ such that
$$\be'(x,y)=\gm(x\trid y)\be(x,y)\gm(y)^{-1}.$$
A non-abelian $2$-cocycle is said to be a \emph{quandle cocycle} if $\be(x,x)=1\ \forall x\in X$.
We denote by $H^2(X,H)$ the set of cohomology classes of $2$-cocycles.
We denote by $H^2_Q(X,H)$ the set of cohomology classes of quandle $2$-cocycles.

If $X$ is a rack, $S$ is a set and $\be:X\times X\to\Ss_S$ is a $2$-cocycle, then there is a
structure of rack in the product $X\times S$ given by
$$(x,s)\trid(y,t)=(x\trid y,\be(x,y)(t)).$$
We denote by $X\times_\be S$ the rack with this structure. Two cohomologous cocycles give rise to
isomorphic structures. If $X$ is a quandle then $X\times S$ is a quandle iff $\be$ is a quandle cocycle.

Let $X$ be a rack, let $\io:X\to X$ be defined by $x\trid \io(x)=x$. Define $(X,\trid^\io)$ as
$X$ with the structure $x\trid^\io y=x\trid\io(y)$. It can be seen that $\io$ is bijective and
$(X,\trid^\io)$ is a quandle.

For a natural number $n$ we denote its $p$-valuation by $v_p(n)$, i.e.,
$v_p(n)=r$ if $n=p^rq$ where $q$ is coprime to $p$.
\end{defin}

\medskip
\section{Racks of order $p^2$}
\begin{theor}\label{et}
Let $X$ be an indecomposable rack of order $p^2$, $p$ a prime number. Then either:
\begin{flalign}
\label{cla1} &X=\ZZ_p\oplus\ZZ_p, &&
	(x_1,x_2)\trid (y_1,y_2)=((1-\al)x_1+\al y_1,(1-\be)x_2+\be y_2) &&
	\al,\be\in\ZZ_p^*-\{1\} \\
\label{cla2} &X=\ZZ_p\oplus\ZZ_p, &&
	(x_1,x_2)\trid (y_1,y_2)=((1-\al)x_1+\al y_1,(1-\al)x_2+\al y_2+y_1-x_1) &&
	\al\in\ZZ_p^*-\{1\} \\
\label{cla3} &X=\FF_{p^2}, &&
	x\trid y=(1-\al)x+\al y &&
	\al\in\FF_{p^2}-\FF_p \\
\label{cla4} &X=\ZZ_{p^2}, &&
	x\trid y=(1-\al)x+\al y &&
	\al\not\equiv 0,1\mod p \\
\label{cla5} &X=\ZZ_{p^2}, &&
	x\trid y=y+1 && \\
\label{cla6} &X=\ZZ_p\oplus\ZZ_p, &&
	(x_1,x_2)\trid(y_1,y_2)=((1-\al)x_1+\al y_1,y_2+1) &&
	\al\in\ZZ_p^*-\{1\}
\end{flalign}
Two racks in different rows are not isomorphic. The non trivial isomorphisms
inside each row are as follows: in \eqref{cla1}, the rack associated to
$(\al,\be)$ is isomorphic to that associated to $(\be,\al)$; in \eqref{cla3}
the rack associated to $\al$ is isomorphic to that associated to $\sg(\al)$,
where $\sg$ is the non-trivial element of the Galois group $\Gal(\FF_{p^2}|\FF_p)$.
\end{theor}
\begin{proof}
If $X$ is faithful then $X$ is a quandle and we prove in \ref{qp2} below that
it is affine. If $X$ is not faithful then we consider its associated quandle
$(X,\trid^{\io})$, which is neither faithful.
Then either $\phi(X)$ has order $1$ or $p$. In the first case $(X,\trid^{\io})$ is
trivial, and then $(X,\trid)$ is given by a permutation $\sg\in\Ss_{p^2}$:
$x\trid y=\sg(y)$. Since $(X,\trid)$ is indecomposable, $\sg$ must be a $p^2$-cycle.
Then it is of the form \eqref{cla5}. In the second case, we have that $(X,\trid^\io)$
is a nonabelian extension of $\phi(X)$ by some set $S$ of order $p$, i.e.,
$(X,\trid^{\io})\simeq (\ZZ_p,\al)\times S$. Since $\phi_x=\phi_{\io x}$,
we have that $\io$ restricts to the fibers, i.e.,
$\io_x=\io|_{\{x\}\times S}:\{x\}\times S\to \{x\}\times S$.
The structure in $X$ can be recovered from $\trid^\io$ as
$(x,s)\trid(y,t)=(x,s)\trid^\io(\io^{-1}(y,t))=((1-\al)x+\al y,\io_y^{-1}(t))$.
Now, it is easy to see that \eqref{Q2} implies
$\io^{-1}_{y\trid z}\io^{-1}_z=\io^{-1}_{x\trid z}\io^{-1}_z$ $\forall x,y,z$,
and thus $\io^{-1}_x=\io^{-1}_y$ $\forall x,y$. We can call then $f=\io^{-1}_x$,
and we have $(x,s)\trid (y,t)=(x\trid y,f(t))$. But for this rack to be
indecomposable $f$ must be a $p$-cycle whence $X$ is isomorphic to the
rack in \eqref{cla6}.

To see that \eqref{cla1}, \eqref{cla2}, \eqref{cla3}, \eqref{cla4} cover
all the affine cases, simply notice that there are two groups of order
$p^2$: $\ZZ_p\oplus\ZZ_p$ and $\ZZ_{p^2}$.
For $\ZZ_p\oplus\ZZ_p$ the isomorphism $g\in\GL_2(\ZZ_p)$ can be either diagonalizable
(class \eqref{cla1}), it can be given by a Jordan block (class \eqref{cla2}) or
its minimal polynomial can be irreducible over $\ZZ_p$ (class \eqref{cla3}).
For $\ZZ_{p^2}$ any automorphism is given by an element in $\ZZ_{p^2}^*$, and we
get class \eqref{cla4}. The conditions on $\al$ and $\be$ in the statement are equivalent
for these quandles to be indecomposable.

Now, by \cite[Lemma 1.33]{ag} two indecomposable affine quandles $(A,g)$ and $(B,h)$
are isomorphic iff there is an isomorphism of the pairs $(A,g)$ and $(B,h)$; i.e, iff
there exists an isomorphism $T:A\to B$ such that $Tg=hT$. This proves that the classes
have no intersection and shows also that the isomorphisms inside each class are
those in the statement.
\end{proof}

Before dealing with the rest of the proof, we derive two corollaries:
\begin{corol}
If $X$ is an indecomposable rack of order $p^2$ then $v_p(|G_X|)=2$ or $3$.
\end{corol}
\begin{proof}
It is proved in \cite{ag} that if $(A,g)$ is an affine quandle then
$G_X^0$ is the semidirect product $A\rtimes (g)$, where $(g)$ is the
subgroup of $\Aut(A)$ generated by $g$. Since the racks in classes
\eqref{cla1}, \eqref{cla2}, \eqref{cla3} and \eqref{cla4} are affine
and the underlying group has order $p^2$, we must prove that
$v_p(|(g)|)\le 1$, but this follows easily by inspection
(actually, this result is part of the proof of \ref{qp2},
which in turn is part of the proof of \ref{et}).

For the rack of class \eqref{cla5} it is clear that $G_X^0$ is cyclic
of order $p^2$, and for those of class \eqref{cla6} it can be seen
that $G_X^0$ is isomorphic to the direct product between a cyclic
group of order $p$ and $G_Y^0$, where $Y$ is the quotient
$(x_1,x_2)\mapsto x_1$. Since $Y\simeq (\ZZ_p,\al)$,
$v_p(|G_Y^0|)=1$ and we are done.
\end{proof}

\begin{corol}
The isomorphisms classes of indecomposable racks of order $p^2$ and their cardinalities
are given in the following table:
\begin{center}
\begin{tabular}{|c|c|c|} \hline
Type & Class & \# \\ \hline
Affine quandle over $\ZZ_p^2$; diagonalizable isomorphism 
	& \eqref{cla1} & $\frac 12(p^2-3p+2)$ \\ \hline
Affine quandle over $\ZZ_p^2$; Jordan block 
	& \eqref{cla2} & $p-2$ \\ \hline
Affine quandle over $\ZZ_p^2$; irreducible polynomial (simple)
	& \eqref{cla3} & $\frac 12p(p-1)$ \\ \hline
Affine quandle over $\ZZ_{p^2}$
	& \eqref{cla4} & $p^2-2p$ \\ \hline
Rack which is not a quandle with $|\phi(X)|=1$
	& \eqref{cla5} & $1$ \\ \hline
Rack which is not a quandle with $|\phi(X)|=p$
	& \eqref{cla6} & $p-2$ \\ \hline
\end{tabular}
\end{center}
\hfill\qed
\end{corol}

We now finish the proof of \ref{et}.
\begin{propo}\label{qp2}
Indecomposable quandles of order $p^2$ are affine. In particular, they are faithful.
\end{propo}
\begin{proof}
Since for $p=2$ this is known, we can assume that $p\neq 2$ (the same tools here
work for the case $p=2$, though sometimes the formulas are easier if we have $\frac 12\in\ZZ_p$).

For simple quandles the result is a consequence of \cite[Thm. 3.12]{ag}. Let $X$ be an indecomposable
non-simple quandle of order $p^2$. Then by \cite[Cor. 2.10]{ag} we have $X\simeq Y\times_{\al}S$,
where $Y$ is an indecomposable quandle of order $p$, $S$ is a set of order $p$ and $\al$ is
a \emph{dynamical $2$-cocycle} (see \cite{ag} for the definition, it is not needed to understand
the proof, though). For $y\in Y$, let us denote by $X_y$ the fibers $y\times S$. These are
quandles and since $Y$ is indecomposable they are all isomorphic, i.e., $X_y\simeq X_{y'}$ as
quandles. We claim that either $X_y$ is indecomposable or either it is trivial. To see this,
take for each $y\in Y$ the decomposition $X_y=\sqcup_nX_y^n$, where $X_y^n$ is the union of
the orbits of $X_y$ with cardinality $n$. Take $(y,s)\in X$; since
$\phi_{(y,s)}:X_{z}\to X_{y\trid z}$ is a quandle isomorphism, it must send $X_z^n$ to
$X_{y\trid z}^n$. Thus, we have a decomposition of $X$ as $X=\sqcup_n(\sqcup_yX_y^n)$.
But $X$ is indecomposable, hence all the orbits in $X_y$ (any $y$) have the same cardinality.
And since $S$ has a prime cardinality, either there is one orbit of order $p$ (and $X_y$ is
indecomposable) or there are $p$ orbits of order $1$ (and $X_y$ is trivial), proving the claim.

We suppose first that $X$ is faithful. By \cite[Thm. A.2]{egs} the group $G_X^0$ is an extension
of a cyclic group by a $p$-group. That is, $G_X^0=N\rtimes_{f} C$, $N$ a $p$-group, $C$ cyclic
and $f$ a $2$-cocycle. Let $M$ be the order of $C$ and let $t$ be a generator of it.
We have the structure $(a,t^i)(b,t^j)=(a\al^i(b)f(i,j),t^{i+j})$,
where $f(i,j)=1\in N$ if $i+j<M$ and $t$ acts by $\al$ on $N$.
Now, the image of $\phi:X\to G_X^0$ is invariant, then it is included in
$\{(a,t^i)\ |\ i=i_0\}$ for some $i_0$. Since this image must generate $G_X^0$, we can assume
that $i_0=1$ (otherwise we re-name $t$ to $t^{i_0}$). The structure of $X$ is given then by
\begin{equation}\label{eq:st}
(a,t)\trid(b,t)=(a,t)(b,t)(a,t)^{-1}=(a,t)(b,t)(1,t^{-1})(a^{-1},1)
	=(a\al(ba^{-1}),t).
\end{equation}
Furthermore, $G_X^0$ has trivial center; in particular if $g\in Z(N)$ then $\al(g)\neq g$,
and, since any $p$ group has a non-trivial center, $\al\neq\id$.
Thus, to classify $X$ we can seek what pairs $N,\al$ can arise and then look for the
structure of the orbits for the action $a\trid b=a\al(ba^{-1})$.
The strategy of the proof in \cite{egs} is the same as this one; however in that case
$X$ can be seen as an orbit in the symmetric group $\Ss_p$ (actually, the faithfulness
condition is immediate), whence the order of $N$ divides the power of $p$ in $\Ss_p$,
which is $p$; thus $N\simeq\ZZ_p$. With a similar reasoning we can prove that the order of $N$
here divides the power of $p$ in $\Ss_{p^2}$, but this is $p^{p+1}$, which gives too much
freedom to $N$. However, the group $N$ is rather small: we claim that $|N|\le p^3$. To see this,
if $X_y$ is indecomposable we have by \cite{egs} that it is affine and
$v_p(|\Aut(X_y)|)=v_p(|\ZZ_p\rtimes\ZZ_p^*|)=1$.
If $X_y$ is not indecomposable, we have that it is trivial and then
$v_p(|\Aut(X_y)|)=v_p(|\Ss_{X_y}|)=1$.
Since $X$ is indecomposable, $Y$ is indecomposable and again $v_p(|G_Y^0|)=1$.
By \cite[Lemma 1.13]{ag} we have a morphism of groups $G^0(\pi):G_X^0\to G_Y^0$ induced by
the projection $\pi:X\to Y$. We look to its kernel $K=\ker G^0(\pi)$. By an abuse of notation,
we denote the elements of $Y$ by those of $\ZZ_p$. Let $w\in K$,
we have $w(X_y)=X_y\ \forall y\in Y$. Then we can restrict $w$ to $X_0$ and $X_1$, i.e.,
we have a homomorphism of groups $R:K\to\Aut(X_0)\times\Aut(X_1)$. But $X_0\cup X_1$ generates
$X$ as a quandle, since $0$ and $1$ generate $Y$. Then $R$ is injective, which proves that
the order of $K$ divides that of $\Aut(X_0)\times\Aut(X_1)=\Aut(S)\times\Aut(S)$.
Then $v_p(|K|)\le 2v_p(|\Aut(S)|)=2$, and $v_p(|G_X|)\le v_p(|K|)+v_p(|G_Y|)\le 3$, proving
the claim.

If $N$ is abelian we are done, since in this case \eqref{eq:st} defines an affine structure.
Thus, if $N$ has order $p$ or $p^2$, there is nothing else to prove. The classification of
groups of order $p^3$ is well known (cf. \cite{b}); there are $3$ abelian groups and two
groups which are not abelian. We must concentrate the attention on the later. We prove in
sections \ref{ss:g1} and \ref{ss:g2} that for each of them we get affine quandles.

\medskip
Suppose now that $X$ is not faithful. The image of $\phi:X\to G_X^0$ must have order $p$, since
if it was $1$ then $X$ would be trivial. Let $Y=\phi(X)$; as before it is an indecomposable quandle
and then $Y\simeq (\ZZ_p,q)$, where $q\in\ZZ_p^*-\{1\}$. By \cite[Prop. 2.20]{ag}, we have
$X=Y\times_{\beta}S$, where $S$ is a set of order $p$ and $\beta:Y\times Y\to\Ss_S$
is a nonabelian quandle $2$-cocycle.
Now, we prove in \ref{h2zp} that this set is trivial; whence any of these quandles is
isomorphic to the product $Y\times S$, $S$ a trivial quandle; but this implies that
$Y\times S$ is decomposable.
\end{proof}

\medskip
\section{The group $(\zpzp)\rtimes C_p$}\label{ss:g1}
Let $G=(\zpzp)\rtimes C_p$ be the group with $t$ a generator of
$C_p$ acting on $\zpzp$ by $t(x,y)t^{-1}=(x,x+y)$.
We denote the elements of $G$ as $(x,y)t^z\ x,y,z\in\ZZ_p$.
Let $\al\in\Aut(G)$. Since $(0,1)$ generates the center of $G$, we must
have $\al(0,1)=(0,q)$ for some $q\in\ZZ_p^*$. Since $\al$ must act
non trivially on the elements of the center, $q\neq 1$.
We denote $\al((1,0))=(a,j)t^b$, $\al(t)=(c,k)t^d$. We have
$\al((x,y))=((a,j)t^b)^x(0,yq)=(xa,xj+ab\frac {x(x-1)}2+yq)t^{bx}$ and
$\al(t^z)=((c,k)t^d)^z=(zc,zk+cd\frac {z(z-1)}2)t^{dz}$, whence
$$\al((x,y)t^z)=(xa+zc,bcxz+xj+zk+ab\frac {x(x-1)}2+cd\frac {z(z-1)}2+yq)t^{bx+dz}.$$

For $\al$ to be a group homomorphism we must have $\al(t)\al((1,0))=\al((1,1))\al(t)$, while
\begin{align*}
\al(t)\al((1,0)) &= (c,k)t^d(a,j)t^b=(c+a,k+j+da)t^{d+b} \\
\al((1,1))\al(t) &= (a,j+q)t^b(c,k)t^d=(c+a,k+j+q+bc)t^{b+d},
\end{align*}
whence necessarily $q=da-bc$.

As said, for $h\in G$, we consider the orbits $\orb_h$ under the action
$g\trid h=g\al(hg^{-1})$. We seek the conditions on $\al$ that give
orbits of order $p^2$. Take $g=1$; we get that
$\{\al^n(h)\ |\ 0\le n<N\}\subseteq\orb_h$. Take $g=(0,y)$ and notice that
$(0,y)\trid h=(0,y)\al(h)(0,-qy)=\al(h)(0,(1-q)y)$. Then
$\{\al^n(h)(0,*)\}\subseteq\orb_h$.
Let us compute the orbit of $\ze=(0,1)$. We have $\{(0,*)\}\subseteq\orb_\ze$.
Acting by $(-x,0)t^{-z}$, we get
\begin{align*}
(-x,0)t^{-z}\trid (0,*) &=(-x,0)t^{-z}\al((0,*)t^z(x,0))=(-x,0)t^{-z}\al((x,*)t^z) \\
	&=(-x,0)t^{-z}(xa+zc,*)t^{bx+dz}=((a-1)x+cz,*)t^{bx+(d-1)z}.
\end{align*}
Let $A=\mdpd acbd$. We get that $\orb_\ze=G$ if $A-1$ is invertible.
Then, $A-1$ is degenerate.

We suppose first that $A=1$. We have $\orb_\ze$ of order $p$. Consider other orbits:
\begin{align*}
(-x,0)t^{-z}\trid(\la,*)t^\mu &= (-x,0)t^{-z}\al((\la,*)t^{\mu+z}(x,0))
	=(-x,0)t^{-z}\al((\la+x,*)t^{\mu+z}) \\
	&=(-x,0)t^{-z}(\la+x,*)t^{\mu+z}=(\la,*)t^{\mu},
\end{align*}
whence all orbits have order $p$. This is not interesting for us.

Suppose then that $A\neq 1$ and that the first row of $A-1$ is non-trivial.
Then $A=\mdpd {r+1}s{ur}{us+1}$ for some $r,s,u\in\ZZ^p$. We compute the
orbit of $(\la,*)$:
\begin{align*}
(-x,0)t^{-z}\trid(\la,*)t^\mu &= (-x,0)t^{-z}\al((\la,*)t^{\mu+z}(x,0))
		=(-x,0)t^{-z}\al((\la+x,*)t^{\mu+z}) \\
	&=(-x,0)t^{-z}((r+1)(\la+x)+s(\mu+z),*)t^{ur(\la+x)+(us+1)(\mu+z)} \\
	&=((\la r+xr+\mu s+zs)+\la,*)t^{u(\la r+xr+\mu s+zs)+\mu},
\end{align*}
whence the orbits are characterized as $\orb^C:=\{(\la,\sigma)t^{\tau}\ |\ \tau-u\la=C\}$,
and this case is interesting for us.

\medskip
If the first row of $A-1$ is trivial we have $A=\mdpd 10r{s+1}$ and the orbits are
\begin{align*}
(-x,0)t^{-z}\trid(\la,*)t^\mu &= (-x,0)t^{-z}\al((\la,*)t^{\mu+z}(x,0))
		=(-x,0)t^{-z}\al((\la+x,*)t^{\mu+z}) \\
	&=(-x,0)t^{-z}((\la+x),*)t^{r(\la+x)+(s+1)(\mu+z)}
		=(\la,*)t^{r\la+rx+s\mu+sz+\mu},
\end{align*}
whence the orbits are characterized as $\orb^C:=\{(\la,\sigma)t^{\tau}\ |\ \la=C\}$,
and this case is also interesting for us.

\medskip
We return then to the case $A=\mdpd{r+1}s{ur}{us+1}$. Let $\Dt_x=\bar x-x$, $\Dt_y=\bar y-y$.
Notice that $q=1+r+us$. We have in $\orb^C$:
\begin{align*}
((x&,y)t^{C+ux})\trid((\bar x,\bar y)t^{C+u\bar x})
	=((x,y)t^{C+ux})\al((\bar x,\bar y)t^{u(\bar x-x)}(-x,-y)) \\
	&=((x,y)t^{C+ux})\al((\bar x-x,\bar y-y-xu(\bar x-x))t^{u(\bar x-x)}) \\
	&=(x,y)t^{C+ux}((\bar x-x)(r+1)+su(\bar x-x),u^2rs(\bar x-x)^2+(\bar x-x)j
		+u(\bar x-x)k+(r+1)ur\frac{(\bar x-x)(\bar x-x-1)}2 \\
	&\hspace*{2cm}	+s(us+1)\frac{u(\bar x-x)(u\bar x-ux-1)}2+(\bar y-y)q-xu(\bar x-x)q)
		\ t^{ur(\bar x-x)+(us+1)u(\bar x-x)} \\
	&=((\bar x-x)(r+1)+su(\bar x-x)+x,u^2rs(\bar x-x)^2+(\bar x-x)j
		+u(\bar x-x)k+(r+1)ur\frac{(\bar x-x)(\bar x-x-1)}2 \\
	&\hspace*{2cm}	+s(us+1)\frac{u(\bar x-x)(u\bar x-ux-1)}2+(\bar y-y)q-xu(\bar x-x)q \\
	&\hspace*{2cm}	+y+(C+ux)((\bar x-x)(r+1)+su(\bar x-x)))
		\ t^{ur(\bar x-x)+(us+1)u(\bar x-x)+C+ux} \\
	&=(q\Dt_x+x\ ,\ \Dt_x^2(u^2rs+\frac u2(r+1)r+\frac{u^2}2s(us+1)) \\
	&\hspace*{1.8cm}+\Dt_x(j+uk-\frac u2(r+1)r-\frac u2s(us+1)+Cq)+q\Dt_y+y)\ t^{uq\Dt_x+C+ux} \\
	&=(q\Dt_x+x\ ,\ \Dt_x^2(\frac u2rq+\frac{u^2}2sq)
	+\Dt_x(j+uk-\frac u2(r+1)r-\frac u2s(us+1)+Cq)+q\Dt_y+y)\ t^{uq\Dt_x+C+ux} \\
\end{align*}

Now we use the following bijection:
$f:\orb^C\to\zpzp$, $f((x,y)t^{ux+C})=(x,y-u\frac{x(x-1)}2)$.
We compute then the quandle structure of $\orb^C$ transported to $\zpzp$
by $f$. We have
\begin{align*}
f(f^{-1}(x,y) &\trid f^{-1}(\bar x,\bar y))
	= f((x,y+u\frac{x(x-1)}2)t^{ux+C}\trid(\bar x,\bar y+u\frac{\bar x(\bar x-1)}2)t^{u\bar x+C}) \\
	&=f((q\Dt_x+x\ ,\ \Dt_x^2(\frac u2rq+\frac{u^2}2sq)
		+\Dt_x(j+uk-\frac u2(r+1)r-\frac u2s(us+1)+Cq)+q\Dt_y+y \\
	&\hspace*{1.8cm}+\frac u2q(\bar x^2-\bar x-x^2+x)+\frac u2x(x-1)) t^{uq\Dt_x+C+ux}) \\
	&=(q\Dt_x+x\ ,\ \Dt_x^2(\frac u2rq+\frac{u^2}2sq+\frac u2q)
		+\Dt_x(j+uk-\frac u2(r+1)r-\frac u2s(us+1)+Cq)+q\Dt_y+y \\
	&\hspace*{1.8cm}+\frac u2q(2\bar xx+2x^2-\bar x+x)+\frac u2x(x-1)-\frac u2(q\Dt_x+x)(q\Dt_x+x-1)) \\
	&=(q\Dt_x+x\ ,\ \Dt_x^2(\frac u2q^2-\frac u2q^2)+q\Dt_y+y \\
	&\hspace*{1.8cm}+\Dt_x(j+uk-\frac u2qr+\frac u2usr-\frac u2sq+\frac u2sr+Cq-\frac u2q+\frac u2q)
		+\frac u2q(2x\Dt_x-2x\Dt_x)) \\
	&=(q\Dt_x+x\ ,\ q\Dt_y+y+\Dt_x(j+uk-\frac u2qr+\frac u2usr-\frac u2sq+\frac u2sr+Cq))
\end{align*}
This means that $\orb^C$ is isomorphic to the affine quandle $(\zpzp,g)$, where
$$g(x,y)=(qx,qy+(j+uk-\frac u2qr+\frac u2usr-\frac u2sq+\frac u2sr+Cq)x).$$

We consider now the case $A=\mdpd 10r{s+1}$. Notice that here $q=s+1$.
Let $\Dt_z=\bar z-z$. We have for the orbit $\orb^C$:
\begin{align*}
(C,y)t^z\trid (C,\bar y)t^{\bar z}
	&=(C,y)t^z\al((C,\bar y)t^{\bar z-z}(-C,\bar y))
		=(C,y)t^z\al((0,\Dt_y-\Dt_zC)t^{\Dt_z}) \\
	&=(C,y)t^z(0,\Dt_zk+q\Dt_y-Cq\Dt_z)t^{(s+1)\Dt_z} \\
	&=(C,q\Dt_y+y+\Dt_z(k-Cq))t^{q\Dt_z+z}, \\
\end{align*}
Then, taking the bijection $f:\orb^C\to\zpzp$ given by $f((C,y)t^z)=(y,z)$ we get
on $\zpzp$ the affine quandle $(\zpzp,g)$ with $g$ given by
$g(y,z)=(qy+(k-Cq)z,qz)$.

\medskip
\section{The group $\ZZ_{p^2}\rtimes C_p$}\label{ss:g2}
Let $G=\ZZ_{p^2}\rtimes C_p$, where $C_p$ is generated by $t$ and the action is given by
$tat^{-1}=a(p+1)$.
Notice that $(at^b)^n=a(n+pb\frac{n(n-1)}2)t^{bn}$.
In particular $(at^b)$ has order $p$ iff $p|a$.

Take $\al\in\Aut(G)$, $\al(1)=at^b$, $\al(t)=ct^d$. Let's compute the conditions
on $a,b,c,d$ for $\al$ to be a homomorphism:
\begin{align*}
\al(t)\al(1)&=(ct^d)(at^b)=(c+(1+pd)a)t^{d+b}=(c+a+pad)t^{d+b} \\
\al(1+p)\al(t)&=(a(1+p))t^b(ct^d)=(a+ap+c+bpc)t^{b+d}
\end{align*}
i.e., we must have $pad=pa+pbc\mod p^2$, or $a(d-1)=bc\mod p$.
On the other hand, $\al(t)$ must have order $p$, whence
$c=0\mod p$. This means that either $a=0\mod p$ or $d=1\mod p$.
The first possibility is excluded since $\al(1)$ must have
order $p^2$. Then, $d=1$ and $\al(t)=(pc')t$.
\begin{align*}
\al(nt^m)&=(at^b)^n(pc't)^m=(a(n+pb\frac{n(n-1)}2))t^{bn}(pc'm)t^m \\
&=(an+pab\frac{n(n-1)}2+pc'm(1+pbn))t^{bn+m} \\
&=(an+p(ab\frac{n(n-1)}2+c'm))t^{bn+m}
\end{align*}
We have, then
\begin{align*}
(nt^m)\trid(\bar nt^{\bar m})
	&=(nt^m)\al(\bar nt^{\Dt_m}(-n))
		=nt^m\al(\bar n-n(1+p\Dt_m)t^{\Dt_m})
		=nt^m\al((\Dt_n-pn\Dt_m)t^{\Dt_m}) \\
	&=nt^m(a(\Dt_n-pn\Dt_m)+p(ab\frac{(\Dt_n-pn\Dt_m)(\Dt_n-pn\Dt_m-1)}2+c'\Dt_m))t^{b\Dt_n+\Dt_m} \\
	&=(n+(1+mp)(a\Dt_n+p(-an\Dt_m+ab\frac{\Dt_n(\Dt_n-1)}2+c'\Dt_m)))t^{b\Dt_n+\bar m} \\
	&=(n+a\Dt_n+p(am\Dt_n+(c'-an)\Dt_m+ab\frac{\Dt_n(\Dt_n-1)}2))t^{b\Dt_n+\bar m} \\
	&=(\bar n-(1-a)\Dt_n+p(c'\Dt_m+a(m\bar n-n\bar m)+ab\frac{\Dt_n(\Dt_n-1)}2))t^{b\Dt_n+\bar m} \\
	&=(\bar n-(1-a)\Dt_n+p(c'\Dt_m+a(-\Dt_m\bar n+\Dt_n\bar m)+ab\frac{\Dt_n(\Dt_n-1)}2))t^{b\Dt_n+\bar m}.
\end{align*}
We see then that if $b\neq 0$ we have
$$\orb_{\bar nt^{\bar m}}\subseteq\{xt^y\ |\ x=\frac{a-1}b(y-\bar m)+a\bar n\mod p\},$$
while if $b=0$ we have
$$\orb_{\bar nt^{\bar m}}\subseteq\{xt^y\ |\ y=\bar m\}.$$
Thus, all orbits have order $\le p^2$. We exclude now the cases in which the order is $<p^2$.
We write $n=r+ps$ ($p\nmid r$), and $D=\bar n-r$, we get
\begin{align*}
(nt^m)\trid(\bar nt^{\bar m})
	&=(\bar n-(1-a)D+p((1-a)s+c'\Dt_m+a(-\Dt_m\bar n+D\bar m)+ab\frac{D(D-1)}2))t^{bD+\bar m} \\
	&=(\bar n-(1-a)D+p((1-a)s+(c'-a\bar n)\Dt_m+aD\bar m+ab\frac{D(D-1)}2))t^{bD+\bar m} \\
\end{align*}
Thus, if $a=1\mod p$, we must have $b\neq 0$ and $c'-\bar n=c'-a\bar n\neq 0\mod p$.

Suppose first that $a\neq 1\mod p$. Then we have to consider the orbits
$\orb^C=\{xt^{-\frac b{1-a}x+C}\}$. We get
\begin{align*}
xt^{\frac {-b}{1-a}x+C}\trid \bar xt^{\frac {-b}{1-a}\bar x+C}
	&=(\bar x-(1-a)\Dt_x+p(c'\frac {-b}{1-a}\Dt_x+a(\frac b{1-a}\Dt_x\bar x-\Dt_x\frac b{1-a}\bar x+\Dt_xC) \\
		&\hspace*{4cm}	+ab\frac{\Dt_x(\Dt_x-1)}2))t^{b\Dt_x-\frac b{1-a}\bar x+C} \\
	&=(\bar x-(1-a)\Dt_x+p((\frac {-bc'}{1-a}+aC-\frac{ab}2)\Dt_x
		+\frac{ab}2\Dt_x^2))t^{b\Dt_x-\frac b{1-a}\bar x+C} \\
\end{align*}
Consider now the function
$$f:\ZZ_{p^2}\to G,\quad f(x)=(x-p\frac b{1-a}\frac{x(x-1)}2)t^{\frac{-b}{1-a}x+C}.$$
We translate the structure of $\orb^C$ to $\ZZ_{p^2}$ via $f$.
We have $f^{-1}(xt^{\frac{-b}{1-a}x+C})=x+p\frac b{1-a}\frac{x(x-1)}2$.
\begin{align*}
f^{-1}(f(x)\trid f(\bar x))
	&=f^{-1}((x-p\frac b{1-a}\frac{x(x-1)}2)t^{\frac{-b}{1-a}x+C} \trid
		(\bar x-p\frac b{1-a}\frac{\bar x(\bar x-1)}2)t^{\frac{-b}{1-a}\bar x+C}) \\
	&=f^{-1}((\bar x-(1-a)\Dt_x+p(-\frac b{1-a}\frac{\bar x(\bar x-1)}2+\frac{b}2(\bar x^2-\bar x-x^2+x) \\
	&\hspace*{2cm}	+(\frac {-bc'}{1-a}+aC-\frac{ab}2)\Dt_x
			+\frac{ab}2\Dt_x^2))t^{b\Dt_x-\frac b{1-a}\bar x+C}) \\
	&=\bar x-(1-a)\Dt_x+p(-\frac b{1-a}\frac{\bar x(\bar x-1)}2+\frac{b}2(\bar x^2-\bar x-x^2+x) \\
	&\hspace*{2cm}	+\frac b{1-a}\frac{(\bar x-(1-a)\Dt_x)(\bar x-(1-a)\Dt_x-1)}2
		+(\frac {-bc'}{1-a}+aC-\frac{ab}2)\Dt_x
			+\frac{ab}2\Dt_x^2) \\
	&=\bar x-(1-a)\Dt_x+p(\frac{b}2(\bar x^2-\bar x-x^2+x)-b\bar x\Dt_x
		+(\frac {-bc'}{1-a}+aC-\frac{ab}2+\frac b2)\Dt_x
			+\frac b2\Dt_x^2) \\
	&=\bar x-(1-a)\Dt_x+p\Dt_x(\frac {-bc'}{1-a}+aC-\frac{ab}2).
\end{align*}
Thus, $\orb^C$ is affine, isomorphic to $(\ZZ_{p^2},g)$,
with $g(x)=(a+p(\frac {-bc'}{1-a}+aC-\frac{ab}2))x$.

Suppose now that $a=1\mod p$. Then we have the orbits $\orb^C=\{(C+px)t^y\}$. We have
$$(C+px)t^y\trid (C+p\bar x)t^{\bar y}=(C+p(x-(1-a)\Dt_x+c'\Dt_y-\Dt_yC))t^{\bar y},$$
but this shows that $\orb^C$ is in this case decomposable as $\sqcup_y\orb^C_y$, with
$\orb^C_y=\{(C+px)t^y\}$.

\medskip
\section{Nonabelian cohomology of the quandle $(\ZZ_p,q)$}
Let $(X,\trid)=(\ZZ_p,q)$ and let $H$ be a group.
Many of the computations here are useful to compute $H^2_Q(X,H)$ when $X=(\FF_{p^t},\omega)$.
\begin{lemma}\label{h2zp}
$H^2_Q(X,H)$ is trivial.
\end{lemma}

\begin{proof}
Let $\be:X\times X\to H$ be a quandle nonabelian $2$-cocycle.
Let $\gm:X\to H$ be defined by
$$\gm(x)=\be(x/(1-q),0)^{-1}.$$
We deform $\be$ by $\gm$, i.e., we take the cohomologous
cocycle $\be'(x,y)=\gm(x\trid y)\be(x,y)\gm(y)^{-1}$. We then have
$$\be'(x,0)=\gm((1-q)x)\be(x,0)\gm(0)^{-1}=\be(x,0)^{-1}\be(x,0)\be(0,0)=1.$$
Also $\be'(x,x)=1\ \forall x\in X$. We then assume that $\be$ has these properties.
The cocycle condition reads as
$$\be((1-q)x+qy,(1-q)x+qz)\be(x,z)=\be(x,(1-q)y+qz)\be(y,z).$$
Take $x=-qz/(1-q)$, and get
$\be(-qz+qy,0)\be(\frac{-qz}{1-q},z)=\be(\frac{-qz}{1-q},(1-q)y+qz)\be(y,z)$,
i.e.,
$$\be(\frac{-qz}{1-q},z)=\be(\frac{-qz}{1-q},(1-q)y+qz)\be(y,z).$$
Take now $y=-qz/(1-q)$, and get
$\be((1-q)x-\frac{q^2z}{1-q},(1-q)x+qz)\be(x,z)=\be(x,0)\be(\frac{-qz}{1-q},z)$,
i.e.,
$$\be((1-q)x-\frac{q^2z}{1-q},(1-q)x+qz)\be(x,z)=\be(\frac{-qz}{1-q},z).$$
In particular,
$\be(\frac{-qz}{1-q},(1-q)x+qz)\be(x,z)=\be((1-q)x-\frac{q^2z}{1-q},(1-q)x+qz)\be(x,z)$,
and then
$$\be(\frac{-qz}{1-q},(1-q)x+qz)=\be((1-q)x-\frac{q^2z}{1-q},(1-q)x+qz).$$
Put now $t=(1-q)x+qz$ and get
$\be(\frac{-qz}{1-q},t)=\be(t-qz-\frac{q^2z}{1-q},t)=\be(t-\frac{qz}{1-q},t)$.
Put $s=-qz/(1-q)$ and get
$$\be(s,t)=\be(t+s,t)\quad\forall s,t\in X.$$
If $t\neq 0$, then $t$ generates $\ZZ_p$ and then $\be(s,t)=\be(t,t)=1$ $\forall s$.
Since for $t=0$ we have $\be(s,t)=1$, we are done.
\end{proof}


\begin{thebibliography}{larguito}
\bibitem [AG]{ag} N. Andruskiewitsch \& M. Gra\~na,
\emph{From racks to pointed Hopf algebras},
\texttt{math.QA/0202084}

\bibitem[B]{b} W. Burnside,
\emph{Theory of groups of finite order},
2nd edition, Dover Pub., New York, 1955.

\bibitem[CJKLS]{cjkls}
J. S. Carter, D. Jelsovsky, S. Kamada, L. Langford and M. Saito,
\emph{State-sum invariants of knotted curves and surfaces from quandle cohomology},
Electron. Res. Announc. Amer. Math. Soc. {\bf 5}
(1999), 146--156 (electronic). Also in {\tt math.GT/9903135}.

\bibitem[EG]{eg} P. Etingof, M. Gra\~na,
\emph{On rack cohomology},
\texttt{math.QA/0201290}

\bibitem[EGS]{egs} P. Etingof, R. Guralnik \& A. Soloviev,
\emph{Indecomposable set-theoretical solutions to the Quantum Yang--Baxter
Equation on a set with prime number of elements},
J. Algebra \textbf{242} (2001), 709-719.

\bibitem[G]{g} M. Gra\~na,
\emph{On Nichols algebras of low dimension}, in
``New Trends in Hopf Algebra Theory";
Contemp. Math. {\bf 267} (2000), pp. 111--136.

\bibitem[J]{j} D. Joyce,
\emph{A Classifying Invariant of Knots, The Knot Quandle},
J. Pure Appl. Alg. \textbf{23} (1982), 37--65.

\bibitem[M]{ma} S. V. Matveev,
\emph{Distributive groupoids in knot theory}
Mat. Sbornik (N.S.) \textbf{119} (161) (1982), no. 1, pp. 78--88, 160. 

\bibitem[Mo]{mo} T. Mochizuki,
\emph{Some calculations of cohomology groups of Alexander quandles},
preprint available at \texttt{http://math01.sci.osaka-cu.ac.jp/~takuro}

\bibitem[O]{os} T. Ohtsuki
\emph{Problems on invariants in knots and $3$-manifolds},
preprint available at \texttt{http://www.is.titech.ac.jp/~tomotada/proj01/problem.ps}


\end{thebibliography}
\end{document}